\documentclass[12 pt]{amsart}
\usepackage{amsmath,times,epsfig,amssymb,amsbsy,amscd,amsfonts,amstext,graphicx,color,bm,amsthm}
\usepackage[all]{xy}
\usepackage{graphicx}

\usepackage[urlcolor=blue]{hyperref}
\usepackage{tikz,subfigure}
\usepackage[vmargin=2.5cm,rmargin=2cm,lmargin=2cm]{geometry}

\usepackage{listings} 

\newlength{\colheight}
\newlength{\colwidth}
\definecolor{main}{rgb}{0.3,0.7,0.9}
\definecolor{red-}{rgb}{1.0, 0.2, 0.0}
\definecolor{blue-}{rgb}{0.0, 0.2, 1.0}
\definecolor{green-}{rgb}{0.0, 0.6, 0.0}
\definecolor{greenblue}{rgb}{0.0, 0.3, 1}
\definecolor{gold}{rgb}{0.8, 0.7, 0.0}

\lstdefinelanguage{cocoa}
{
  basicstyle=\small\ttfamily,
  commentstyle=\color{red!90!black},       
  stringstyle=\color{green!70!black},      
  morecomment=[l]{//},
  morecomment=[l]{--},
  morecomment=[s]{/*}{*/},
  morestring=[b]{"}, 
  classoffset=1,
  morekeywords={
    define,enddefine,if,endif,for,endfor, 
    use,in,then,else,elif,return,
    and,or,
    break,continue,ciao,do,exit,
    ImportByValue,ImportByRef,importbyvalue,importbyref,
    in,isin,IsIn,
    on,opt,PrintLn,println,print,
    quit,ref,return,step,
    toplevel,TopLevel,
    then,to,
    use
  },
  keywordstyle=\color{blue!70!green!80!white},
  classoffset=2,
  morekeywords={
    Ideal,Mat,       Not,Record,Error,
    ideal,mat,matrix,not,record,error,submodule
  },
  keywordstyle=\color{purple!50!white},
  classoffset=3,
  morekeywords={
    TRUE,FALSE,True,False,true,false,
    Lex,Xel,DegLex,DegRevLex,
    PosTo,ToPos,Null
  },
  keywordstyle=\color{brown},
}

\def \X(#1){\{x_1,\dots, x_{#1}\}}

\begin{document}

\title{Enumerative combinatorics package for CoCoA}

\begin{abstract}
We introduce the package \textbf{combinatorics} for the software CoCoA. This package provides a data structure and the necessary methods for computing several known enumerative combinatorial invariants. 
\end{abstract}

\author{Akin Scott}
\address{Akin Scott, Department of Mathematics \& Statistics, Northern Arizona University,
801 S Osborne Drive,
Flagstaff, AZ 86011, USA.}
\email{sa2638@nau.edu}
\author{Michele Torielli}
\address{Michele Torielli, Department of Mathematics \& Statistics, Northern Arizona University,
801 S Osborne Drive,
Flagstaff, AZ 86011, USA.}
\email{michele.torielli@nau.edu}


\date{\today}
\maketitle


\section{Introduction}
Enumerative combinatorics is a branch of combinatorics that focuses on counting the number of ways in which certain patterns, configurations, or discrete structures can be formed. Classic and well-known examples include counting permutations, combinations, and arrangements of objects under various constraints. More generally, given an infinite family of finite sets $\{A_n\}_{n\in\mathbb{N}}$, enumerative combinatorics seeks to describe a counting function that assigns to each natural number $n$ the number of elements in the corresponding set $A_n$. This function captures how the size of the sets evolves as the parameter $n$ changes.
Although determining the size of a set is a very broad mathematical problem, many questions arising in science, engineering, and everyday applications can be translated into relatively simple combinatorial models. These models often allow complex real-world situations to be studied through discrete structures that are amenable to precise counting techniques. In favorable cases, the resulting counting functions can be expressed in the form of closed formulas, that is, explicit expressions constructed from elementary operations such as addition, multiplication, exponentiation, and factorials.
The process of discovering such explicit expressions is known as algebraic enumeration. It typically involves first identifying a recurrence relation that relates the sizes of successive sets, or constructing a generating function that encodes the counting sequence in a compact algebraic form. These intermediate tools can then be analyzed and manipulated to derive a closed-form formula, providing both an exact description of the counting sequence and deeper insight into the underlying combinatorial structure. We refer the
reader to \cite{Sagan} and \cite{Stanley} for a comprehensive account of this subject.

%

In this paper, we describe the new package \textbf{combinatorics} for the software CoCoA (\cite{CoCoALib}, \cite{AbbottBigatti2016} and \cite{COCOA}).
This package computes several known enumerative combinatorial invariants that were previously missing from the CoCoA's library. It makes use of several already implemented functions of CoCoA, like binomial coefficients and permutations, and it complements nicely the already implemented package on posets (see \cite{ArrPackage}).

We introduce the package \textbf{combinatorics} via several examples.  
Specifically, in Section 2 we describe the implemented functions related to the notions of permutations and combinations. In Section 3, we describe all the generalizations of the Catalan numbers that are present in our package. In Section 4, we talk about special types of prime numbers. In Section 5, we describe other functions that compute variouos enumerative invariants, like falling and rising factorial, Bell numbers and Touchard polynomials.

This package will be part of the next official CoCoA release, and it was initiated by the first author as requirement to pass the graduate course MAT 690 at NAU in Spring 2025.

\section{Permutations and similar numbers}

In combinatorics, there are several fundamental structures that are at the core of most counting arguments. A \textbf{permutation} of $n$ indices is a bijective map from the set $[n]:=\{1,\dots, n\}$ to itself. Notice that the set of all permutations $S_n$ on $n$ indices form a group under composition called the permutation or symmetric group. One can then generalize this concept to the notion of a \textbf{$k$-permutation} ( or a \textbf{partial permutation}) on $n$ indices, which is an ordered selection of $k$ distinct items from a set of $n$ distinct indices. A \textbf{$k$-combination} on $n$ indices is a selection of $k$ distinct items from a set of $n$ distinct indices, such that the order of the selection does not matter (unlike $k$-permutations). A \textbf{$k$-sequence} on $n$ indices is an ordered selection of $k$ items, in which repetition are allowed, from a set of $n$ distinct indices. A \textbf{$k$-collection} on $n$ indices is a selection of $k$ items, in which repetition are allowed, from a set of $n$ distinct indices, such that the order of the selection does not matter (unlike $k$-sequences). A \textbf{derangement} on $n$ indices is a permutation on $n$ indices such that no index is mapped to itself.
Because of their central role in the theory, we implemented their construction and counting functions.

\begin{lstlisting}[alsolanguage=cocoa]
/**/ NumPartialPermutations(5,3);
60
/**/ NumCombinations(5,3);
10
/**/ NumCollections(5,3);
35
/**/ NumSequences(5,3);
125
/**/ Derangements(4);
[[2,  1,  4,  3],  [2,  3,  4,  1],  [2,  4,  1,  3],  [3,  1,  4,  2],  [3,  4,  1,  2],  [3,  4,  2,  1],  [4,  1,  2,  3],  [4,  3,  1,  2],  [4,  3,  2,  1]]
/**/ NumDerangements(4);
9
/**/ NumDerangements(6)=len(Derangements(6));
true
\end{lstlisting}

Between all possible permutations, it is important to be able to compute the ones that avoid certain patterns, like permutations that avoid the pattern 132. These are the permutation $\sigma$ on $n$ indices with the property that do not exist $1\le i<j<k\le n$ such that $\sigma(i)<\sigma(k)<\sigma(j)$. These type of permutations have connections with other combinatorial objects. For example the number of the permutations on $n$ indices that avoids 132 coincides with the $n$-th Catalan numbers $C_n$, see \cite{Separable}. Moreover, we also want to be able to compute the \textbf{separable permutations}, i.e., the permutations that can be obtained from the trivial permutation by direct sums and skew sums,  see \cite{Separable}. Moreover, separable permutations may be characterized by the forbidden permutation patterns 2413 and 3142.

\begin{lstlisting}[alsolanguage=cocoa]
/**/ AvoidingPermutations(4,"132");
[[1,  2,  3,  4],  [2,  1,  3,  4],  [2,  3,  1,  4],  [2,  3,  4,  1],  [3,  1,  2,  4],  [3,  2,  1,  4],  [3,  2,  4,  1],  [3,  4,  1,  2],  [3,  4,  2,  1],  [4,  1,  2,  3],  [4,  2,  1,  3],  [4,  2,  3,  1],  [4,  3,  1,  2],  [4,  3,  2,  1]]
/**/ AvoidingPermutations(3,"231");
[[1,  2,  3],  [1,  3,  2],  [2,  1,  3],  [3,  1,  2],  [3,  2,  1]]
/**/ CatalanNumber(7)=len(AvoidingPermutations(7,"132"));
true
/**/ CatalanNumber(9);
4862
/**/ 4862=len(AvoidingPermutations(9,"132"));
true
/**/ SeparablePermutations(4);
[[1,  2,  3,  4],  [1,  2,  4,  3],  [1,  3,  2,  4],  [1,  3,  4,  2],  [1,  4,  2,  3],  [1,  4,  3,  2],  [2,  1,  3,  4],  [2,  1,  4,  3],  [2,  3,  1,  4],  [2,  3,  4,  1],  [2,  4,  3,  1],  [3,  1,  2,  4],  [3,  2,  1,  4],  [3,  2,  4,  1],  [3,  4,  1,  2],  [3,  4,  2,  1],  [4,  1,  2,  3],  [4,  1,  3,  2],  [4,  2,  1,  3],  [4,  2,  3,  1],  [4,  3,  1,  2],  [4,  3,  2,  1]]
\end{lstlisting}

\section{From the Catalan numbers}

Similarly to the Catalan numbers, there are several other numbers that compute important quantities in combinatorics. Since the Catalan numbers can be interpreted as a special case of the Bertrand's ballot theorem, i.e., $C_n$ can be seen as the number of ways for a candidate $A$ with $n+1$ votes to lead candidate $B$ with $n$ votes, see \cite{StanleyCatalan}. This allow us to generalize the notion of Catalan numbers and define $T(n,m)=\frac{(2m)!(2n)!}{(m+n!m!n!)}$. Following Ira Gessel, these numbers are called the \textbf{super Catalan numbers} and they satisfy $T(1,n)=2C_n$. There are other known combinatorial descriptions for $m=2,3$ and $4$, however it is an open problem to find a general combinatorial interpretation for these numbers, see \cite{Chen, Gheo}.
The super Catalan numbers should not confused with the \textbf{ Schr\"oder–Hipparchus numbers}, which sometimes are also called super-Catalan numbers or little Schr\"oder numbers. These numbers can be used to count the plane trees with a given set of leaves, the ways of inserting parentheses into a sequence, and the ways of dissecting a convex polygon into smaller polygons by inserting diagonals, and they are computed using the formula $x_n= \frac{1}{n}\sum_{k=0}^{n-1}2^k\binom{n}{k}\binom{n}{k+1}$, see \cite{Stanley}. 
The \textbf{Fuss-Catalan numbers}, also known as Raney numbers, $A_m(p,r)=\frac{r}{mp+r}\binom{mp+r}{m}$ have many combinatorial interpretations, see \cite{PatternFussCat, Rusu, Tsujie2}. For example they count the number of legal permutations or allowed ways of arranging a number of articles, that is restricted in some way. The \textbf{Catalan-Mersenne numbers} $c_n$ are defined by the recursive formula $c_0=2$ and $c_{n+1}=2^{c_n}-1$. Moreover, these numbers are conjectured by Catalan himself to be prime ``up to a certain limit``. Although the first five terms are prime, no known methods can prove that any further terms are prime (in any reasonable time) simply because they are too huge. The \textbf{Mersenne numbers} $M_n$ are the numbers of the form $2^n-1$, and this numbers can be prime only if $n$ is a prime. The \textbf{double Mersenne numbers} are numbers of the form $M_{M_p}=2^{2^p}-1$ when $p$ is prime, and it is prime only if $M_p$ is already prime. The \textbf{Schr\"oder numbers} $S_n=\frac{1}{n}\sum_{k=0}^n2^k\binom{n}{k}\binom{n}{k-1}$  also called a large Schr\"oder number or big Schr\"oder number, describes the number of lattice paths from the southwest corner $(0,0)$ of a $n\times n$ grid to the northeast corner $(n,n)$, using only single steps north $(0,1)$, northeast $(1,1)$ or east $(1,0)$ that do not rise above the SW-NE diagonal, see \cite{EncIntegerSeq}. The \textbf{Narayana numbers} $N(n,k)=\frac{1}{n}\binom{n}{k}\binom{n}{k-1}$ refine Catalan numbers, in fact $\sum_{k=1}^nN(n,k)=C_n$, and count various objects like Dyck paths with specific peaks or permutations with certain descent patterns, see \cite{Narayama}. Associated to these numbers, we can define the \textbf{Narayana polynomials} defined by $N_n(t)=\sum_{i=0}^nN(n,i)t^i$. In this case $N_n(1)=C_n$ and $N_n(2)=S_n$.

\begin{lstlisting}[alsolanguage=cocoa]
/**/ SuperCatalanNumber(5,3);
90
/**/ SuperCatalanNumber(1,4)=2*CatalanNumber(4);
true
/**/ SchroderHipparchusNumber(3);
11
/**/ SchroderHipparchusNumber(6);
903
/**/ FussCatalanNumber(2,3,3);
12
/**/ CatalanMersenneNumber(4);
170141183460469231731687303715884105727
/**/ IsPrime(CatalanMersenneNumber(4));
true
/**/ MersenneNumber(6);
63
/**/ MersenneNumber(5);
31
/**/ IsPrime(MersenneNumber(3));
true
/**/ DoubleMersenneNumber(3);
127
/**/ DoubleMersenneNumber(5);
2147483647
/**/ 
/**/ SchroderNumber(3);
22
/**/ NarayanaNumber(4,2);
6
/**/ NarayanaNumber(4,4);
1
/**/ NarayanaPoly(4);
t^4 +6*t^3 +6*t^2 +t
/**/  NarayanaPoly(8);
t^8 +28*t^7 +196*t^6 +490*t^5 +490*t^4 +196*t^3 +28*t^2 +t
/**/ CatalanNumber(5)=eval(NarayanaPoly(5),[1]);
true
/**/ SchroderNumber(5)=eval(NarayanaPoly(5),[2]);
true
\end{lstlisting}

There are other numbers that generalize the Catalan ones. In particular, the \textbf{Catalan's triangle} has entries $C(n,k)=\binom{n+k}{k}-\binom{n+k}{k-1}$ and we clearly have that $C(n,n)=C_n$. Moreover, $C(n,k)$ gives the number of strings consisting of $n$ X's and $k$ Y's such that no initial segment of the string has more Y's than X's, see \cite{Shapiro}. 
The \textbf{Catalan's trapezoids} are a set of number trapezoids which generalize the Catalan’s triangle. Catalan's trapezoid of order $m$ is defined as
\begin{align*}
C_m(n,k) &=
  \begin{cases}
 \binom{n+k}{k} & \text{ if $0\le k<m;$} \\
 \binom{n+k}{k}-\binom{n+k}{k-m} & \text{ if $m\le k\le n+m-1;$} \\
 0 & \text{ if $k>n+m-1.$}
  \end{cases}
\end{align*}
They count the number of strings consisting of $n$ X's and $k$ Y's such that in every initial segment of the string the number of Y's does not exceed the number of X's by $m$ or more, see \cite{Shlomi}. By definition, we have that $C_1(n,k)=C(n,k)$, i.e., the Catalan's trapezoid of order $m = 1$ is Catalan's triangle.

\begin{lstlisting}[alsolanguage=cocoa]
/**/ CatalanTriangleNumber(4,4);
14
/**/ CatalanTriangleNumber(5,4);
42
/**/ CatalanTriangleNumber(25,14);
-- 6962078952
/**/ CatalanTrapezoidNumber(3,3,3);
19
/**/ CatalanTrapezoidNumber(7,4,3);
319
/**/ CatalanTrapezoidNumber(37,14,13);
1292706174849
/**/ CatalanTrapezoidNumber(7,4,1)=CatalanTriangleNumber(7,4);
true
/**/ CatalanTrapezoidNumber(15,10,1)=CatalanTriangleNumber(15,10);
true
\end{lstlisting}

\section{Prime numbers}

Among all prime numbers, there are several with additional properties. A prime number $p$ is called a \textbf{factorial prime} if there exists $n$ such that $p=n!\pm 1$. Given an integer $n$ the \textbf{primorial} of $n$ is the product of all primes smaller or equal to $n$. This allows to define the notion of a \textbf{primorial prime}, that is a prime number that is exactly one greater or less than a primorial number.

\begin{lstlisting}[alsolanguage=cocoa]
/**/ FactorialPrime(2); -- the n-th factorial prime
3
/**/ FactorialPrime(5);
23
/**/ FactorialPrimesList(10); -- first n factorial primes
[2,  3,  5,  7,  23,  719,  5039,  39916801,  479001599,  87178291199]
/**/ FactorialPrimesFinder(10); -- factorial primes smaller to n!+1
[2,  3,  5,  7,  23,  719,  5039]
/**/ PrimorialPrime(3); -- n-th primorial prime
5
/**/ PrimorialPrimeList(5); -- first n primorial primes
[2,  3,  5,  7,  29]
/**/ primorial(5); -- already existing CoCoA function
30
/**/ PrimorialPrimesFinder(5); -- primorial primes smaller to 1+ primorial(n)
[2,  3,  5,  7,  29,  31]
/**/ 
\end{lstlisting}

\section{Other numbers}

The \textbf{rising factorial} (sometimes called the Pochhammer function, Pochhammer polynomial, ascending factorial, rising sequential product, or upper factorial) is defined as the polynomial $(x)^n=\prod_{k=0}^{n-1}(x+k)$. The \textbf{falling factorial} (sometimes also called the descending factorial, falling sequential product, or lower factorial) is defined as the polynomial $(x)_n=\prod_{k=0}^{n-1}(x-k)$. It turns out that the coefficients that appear in the expansion of $(x)_n$ are \textbf{Stirling numbers of the first kind} $s(n,k)$. We can also define the \textbf{Stirling number of the second kind} (or Stirling partition number) as the number of ways to partition a set of $n$ objects into $k$ non-empty subsets and is denoted by $S(n,k)$, see \cite{ConcreteMaths}. In addition, we have that $S(n,k)=\sum_{j=0}^k\frac{(-1)^{k-j}j^n}{(k-j)!j!}$. The \textbf{Bell numbers} $B_n=\sum_{k=0}^nS(n,k)$ count the different ways to partition a set that has exactly $n$ elements, or equivalently, the equivalence relations on it, see \cite{Bellnumb}, and they also count the different rhyme schemes for $n$-line poems. Analogously, the \textbf{ordered Bell numbers} or Fubini numbers count the weak orderings on a set of $n$ elements, weak orderings arrange their elements into a sequence allowing ties and they can also be computed from the Stirling numbers of the second kind via $a_n=\sum_{k=0}^nk!S(n,k)$. From the Stirling numbers of the second kind, we can also compute the \textbf{Touchard polynomials}, also called the exponential polynomials, using the formula $T_n(t)=\sum_{k=0}^nS(n,k)t^k$. Notice that $T_n(1)=B_n$.

\begin{lstlisting}[alsolanguage=cocoa]
/**/ use QQ[x];
/**/ RisingFactorial(x,5);
x^5 +10*x^4 +35*x^3 +50*x^2 +24*x
/**/ FallingFactorial(x,5);
x^5 -10*x^4 +35*x^3 -50*x^2 +24*x
/**/  [StirlingNumber(5,k) | k in 1..5];
[24,  -50,  35,  -10,  1]
/**/ StirlingNumber2(4,3);
6
/**/ BellNumber(3);
5
/**/ BellNumber(4);
15
/**/ OrderedBellNumber(3);
13
/**/ OrderedBellNumber(4);
75
/**/ TouchardPoly(3);
t^3 +3*t^2 +t
/**/ TouchardPoly(5);
t^5 +10*t^4 +25*t^3 +15*t^2 +t
/**/ BellNumber(8)=eval(TouchardPoly(8),[1]);
true
\end{lstlisting}

Another fundamental object is the \textbf{Pascal's triangle}, that is an infinite triangular array of the binomial coefficients which play a crucial role in probability theory, combinatorics, and algebra.

\begin{lstlisting}[alsolanguage=cocoa]
/**/ PrintPascalTriangle(5); -- first 5 rows
[1]
[1,  1]
[1,  2,  1]
[1,  3,  3,  1]
[1,  4,  6,  4,  1]
[1,  5,  10,  10,  5,  1]
/**/ PrintPascalTriangle(3,5); -- from row 3 to row 5
[1,  3,  3,  1]
[1,  4,  6,  4,  1]
[1,  5,  10,  10,  5,  1]
\end{lstlisting}

\bigskip
\paragraph{\textbf{Acknowledgements}} 
During the preparation of this paper, the second author was partially supported by the Perko Research Award for 2025-2026. The authors would like to thank Anna Bigatti and John Abbott for helping make the package more efficient.

\bibliographystyle{plain}

\end{document}